\documentclass[12pt]{article}

\usepackage{a4,graphicx,amsmath,amsfonts,amssymb,mathrsfs,bbm,color}

\newcommand{\R}{\mathbbm{R}}

\newcommand{\N}{\mathbbm{N}}
\newcommand{\ltworandom}{\mathcal{L}^2(\mathcal{M},\rho)}
\newcommand{\htwonorm}{\mathcal{H}_2}

\newtheorem{definition}{Definition}
\newtheorem{theorem}{Theorem}

\newtheorem{lemma}{Lemma}

\newtheorem{remark}{Remark}

\begin{document}
\begin{center}
  {\bf \Large Stochastic Galerkin method \\[0.1ex] 
      and port-Hamiltonian form \\[0.7ex]
      for linear dynamical systems of second order}

\vspace{10mm}

{\large Roland~Pulch} \\[1ex]
{\small Institute of Mathematics and Computer Science, 
Universit\"at Greifswald, \\
Walther-Rathenau-Str.~47, 17489 Greifswald, Germany. \\
Email: {\tt roland.pulch@uni-greifswald.de}}

\end{center}

\bigskip\bigskip


\begin{center}
{\bf Abstract}

\begin{tabular}{p{13.5cm}}
We investigate linear dynamical systems of second order. 
Uncertainty quantification is applied, where physical parameters 
are substituted by random variables. 
A stochastic Galerkin method yields a linear dynamical system 
of second order with high dimensionality. 
A structure-preserv\-ing model order reduction (MOR) produces a small 
linear dynamical system of second order again. 
We arrange an associated port-Hamiltonian (pH) formulation of first order 
for the second-order systems. 
Each pH system implies a Hamiltonian function describing an internal energy. 
We examine the properties of the Hamiltonian function for 
the stochastic Galerkin systems. 
We show numerical results using a test example, 
where both the stochastic Galerkin method and structure-preserving MOR 
are applied.  
 
\medskip

\end{tabular}
\end{center}

\medskip


\section{Introduction}
Port-Hamiltonian (pH) formulations of dynamical systems represent an 
advantageous and useful structure, 
see~\cite{jacob-zwart,schaft-jeltsema}. 
There is a close connection of pH systems to passivity and stability, 
see~\cite{beattie-mehrmann-dooren,jacob-skrepek}.
Models derived from pH formulations often inherit 
the beneficial properties. 
Furthermore, each pH system exhibits a Hamiltonian function, 
which characterises an internal energy. 

We consider linear systems of ordinary differential equations (ODEs) 
of second order. 
Such systems are generated by a modelling of mechanical 
mass-spring-damper systems and some types of electric circuits,  
for example, see~\cite{braun,inman}. 
In the linear systems, the coefficient matrices are symmetric and 
positive definite or semi-definite. 
We replace physical parameters of a system by independent random variables 
to perform an uncertainty quantification (UQ). 
A stochastic Galerkin approach, see~\cite{sullivan,xiu-book}, 
generates a larger linear system of ODEs of second order. 
The larger matrices inherit the definiteness of the original matrices, 
i.e., the stochastic Galerkin projection is structure-preserving. 

Both the linear second-order ODEs and their stochastic Galerkin systems 
feature an associated pH formulation of first order. 
We investigate the meaning of the Hamiltonian function in the 
pH system for the stochastic Galerkin approach. 
The relations to the input-output behaviour of both the original 
linear dynamical system and the linear stochastic Galerkin system 
are addressed. 

Furthermore, the linear stochastic Galerkin system represents an 
adequate candidate for a model order reduction (MOR) due to the 
high dimensionality. 
We apply projection-based MOR methods, see~\cite{antoulas}.
Stability-preserving MOR of first-order linear dynamical systems 
was examined in~\cite{castane-selga,ionescu,pulch-jmi,pulch21}, 
for example. 
In MOR of second-order linear dynamical systems, 
structure-preser\-va\-tion is desired, 
where the properties of the reduced matrices are identical, see
\cite{bai-su,beattie-gugercin,chu-tsai-lai,reis-stykel,salimbahrami-lohmann}. 
We achieve stability-preservation by structure-preservation 
in the second-order systems. 
Consequently, each reduced-order model owns an associated pH formulation 
of first order again.   

The article is organised as follows. 
We specify pH systems of first order and the investigated systems 
of ODEs of second order in Section~\ref{sec:modelling}. 
The stochastic modelling as well as the stochastic Galerkin method are 
applied in Section~\ref{sec:stochastic-model}. 
We examine the pH formulation of the stochastic Galerkin projection and 
its Hamiltonian function. 
Section~\ref{sec:mor} includes a structure-preserving MOR of the 
stochastic Galerkin system. 
Finally, we present results of numerical computations for a 
test example in Section~\ref{sec:example}.


\section{Problem Definition}
\label{sec:modelling}
We review port-Hamiltonian systems and the class of linear dynamical systems 
under investigation.

\subsection{Port-Hamiltonian Systems}
$\mathcal{S}_{\succ}^n$ and $\mathcal{S}_{\succeq}^n$ denote
the set of real symmetric square matrices of dimension~$n$, 
which are positive definite and positive semi-definite, respectively.
$I_n \in \R^{n \times n}$ represents the identity matrix. 

We define linear port Hamiltonian systems in the case of implicit ODEs, 
cf.~\cite{beattie-mehrmann}.

\begin{definition} \label{def:ph}
The form of a first-order linear port-Hamiltonian (pH) system is
\begin{align} \label{eq:ph}
\begin{split}
E \dot{x} & = (J-R) Q x + (B-P) u \\
y & = (B+P)^\top Q x + (S+N) u 
\end{split}
\end{align}
with matrices $J,R,Q,E \in \R^{n \times n}$, 
$B,P \in \R^{n \times m}$, $S,N \in \R^{m \times m}$ 
such that~$J$ and~$N$ are skew-symmetric, 
$E^\top Q \in \mathcal{S}_{\succ}^n$, and 
\begin{equation} \label{eq:matrix-w} 
W = \begin{pmatrix} 
Q^\top R Q & Q^\top P \\
P^\top Q & S \\
\end{pmatrix}
\in \mathcal{S}_{\succeq}^{n+m} . 
\end{equation}
The associated Hamiltonian function reads as 
\begin{equation} \label{eq:hamiltonian}
H(x) = \tfrac{1}{2} x^\top (E^\top Q) x 
\end{equation}
for $x \in \R^n$.
\end{definition}

\begin{remark} \mbox{}  
\begin{itemize} \em
\item[i)]
In the case of ODEs, the mass matrix~$E$ is non-singular. 
In the case of dif\-fer\-en\-tial-algebraic equations (DAEs), 
the mass matrix~$E$ is singular and the condition 
$E^\top Q \in \mathcal{S}_{\succ}^n$ has to be weakened into 
$E^\top Q \in \mathcal{S}_{\succeq}^n$, 
see~\cite{hauschild}.
\item[ii)] 
The property~(\ref{eq:matrix-w}) implies that 
$Q^\top R Q \in \mathcal{S}_{\succeq}^n$ and $S\in \mathcal{S}_{\succeq}^n$ 
is symmetric as well.
\end{itemize}
\end{remark}

Note that it is necessary for a pH system to have the 
same number of inputs in~$u$ as the number of outputs in~$y$.
The Hamiltonian function of a pH system represents an internal energy. 
Considering a transient solution of the pH system, it holds that
\begin{align} \label{eq:hamiltonian-bound} 
\begin{split}
H(x(t_2)) - H(x(t_1)) & = 
\int_{t_1}^{t_2} y(t)^\top u(t) - x(t)^\top Q^\top R Q x(t) \; 
\mathrm{d}t \\[1ex]
& \le \int_{t_1}^{t_2} y(t)^\top u(t) \; \mathrm{d}t 
\end{split}
\end{align}
for $t_2 \ge t_1 \ge 0$. 
The latter relation is called the dissipation inequality. 
The inner product $y^\top u$ has units of power, 
representing the rate at which energy is delivered to the system.

\subsection{Linear Dynamical Systems of Second Order}
Let parameters $\mu \in \mathcal{M}$ be given in a domain $\mathcal{M} \subseteq \R^q$.
We consider linear dynamical systems of second order in the form
\begin{align} \label{eq:system-secondorder} 
\begin{split}
M(\mu) \ddot{p} + D(\mu) \dot{p} + K(\mu) p & = B'(\mu) u \\[1ex]
F(\mu) p + G(\mu) \dot{p} & = y
\end{split}
\end{align}
with matrices $M,D,K : \mathcal{M} \rightarrow \R^{n \times n}$, 
$B' : \mathcal{M} \rightarrow \R^{n \times n_{\rm in}}$, and $F,G : \mathcal{M} \rightarrow \R^{n_{\rm out} \times n}$.
Therein, we assume $M(\mu),K(\mu) \in \mathcal{S}_{\succ}^n$ and 
$D(\mu) \in \mathcal{S}_{\succeq}^n$ for all $\mu \in \mathcal{M}$. 
These assumptions are typically satisfied for mass-spring-damper systems, for example. 
If it holds that $D(\mu) \in \mathcal{S}_{\succ}^n$, 
then the linear dynamical system~(\ref{eq:system-secondorder}) 
is asymptotically stable, 
see~\cite{inman}. 
In the case of $D(\mu) \in \mathcal{S}_{\succeq}^n$, the system may be stable or unstable. 
In~(\ref{eq:system-secondorder}), the quantity of interest (QoI) is the 
output $y : [0,\infty) \times \mathcal{M} \rightarrow \R^{n_{\rm out}}$.

Following~\cite{beattie-mehrmann}, 
a system~(\ref{eq:system-secondorder}) can be written in the 
pH form~(\ref{eq:ph}) using the state variables 
$x=(\dot{p}^\top,p^\top)^\top$ and the matrices 
$$ E(\mu) = 
\begin{pmatrix} M(\mu) & 0 \\ 0 & I_{n} \\ \end{pmatrix} , \qquad 
J = 
\begin{pmatrix} 0 & -I_{n} \\ I_{n} & 0 \\ \end{pmatrix} , \qquad 
R(\mu) = 
\begin{pmatrix} D(\mu) & 0 \\ 0 & 0 \\ \end{pmatrix} $$
$$ Q(\mu) = 
\begin{pmatrix} I_{n} & 0 \\ 0 & K(\mu) \\ \end{pmatrix} , \qquad
B(\mu) = 
\begin{pmatrix} B'(\mu) \\ 0 \\ \end{pmatrix} , $$
and $P=0$, $S=0$, $N=0$. 
To attain pH form, the involved output matrices become $F = 0$ and 
$G(\mu) = B'(\mu)^\top$. 
The Hamiltonian function~(\ref{eq:hamiltonian}) is
\begin{equation} \label{eq:hamiltonian-deterministic}
H(x(t,\mu),\mu) = \tfrac{1}{2} \left(
\dot{p}(t,\mu)^\top M(\mu) \dot{p}(t,\mu) + 
p(t,\mu)^\top K(\mu) p(t,\mu) \right)
\end{equation}
for each $\mu \in \mathcal{M}$.

\section{Stochastic Model} 
\label{sec:stochastic-model}
We use a stochastic model for the quantification of uncertainties 
in dynamical systems. 

\subsection{Random Variables and Function Spaces}
The physical parameters~$\mu \in \mathcal{M} \subseteq \R^q$, 
which are included in a system~(\ref{eq:system-secondorder}), 
are replaced by independent random variables on a probability space 
$(\Omega,\mathcal{A},P)$. 
Traditional probability distributions can be employed for each 
random variable like uniform distribution, beta distribution, 
Gaussian distribution, etc. 
We assume that a joint probability density function $\rho : \mathcal{M} \rightarrow \R$ 
is available. 
The expected value of a measurable function $f : \mathcal{M} \rightarrow \R$ reads as 
\begin{equation} \label{eq:expected-value}
\mathbb{E}[f] = \int_{\Omega} f(\mu(\omega)) \; \mathrm{d}P(\omega)  
= \int_{\mathcal{M}} f(\mu) \rho(\mu) \; \mathrm{d}\mu .
\end{equation}
Likewise, we obtain an inner product for two measurable functions $f,g$
\begin{equation} \label{eq:inner-product}
\langle f , g \rangle = 
\int_{\mathcal{M}} f(\mu) g(\mu) \rho(\mu) \; \mathrm{d}\mu . 
\end{equation} 
This inner product is defined on the space of square integrable functions 
$$ \ltworandom = \left\{ f : \mathcal{M} \rightarrow \R \; : \; 
f \; \mbox{measurable and} \; \mathbb{E}[f^2] < \infty \right\} . $$
We denote the dedicated norm by $\| f \|_{\ltworandom} = \sqrt{\langle f , f \rangle}$. 

\subsection{Polynomial Chaos Expansions}
Let $(\Phi_i)_{i \in \N}$ be an orthonormal basis with respect to the 
inner product~(\ref{eq:inner-product}), which consists of multivariate
polynomials $\Phi_i : \mathcal{M} \rightarrow \R$. 
We assume that $\Phi_1 \equiv 1$ is the unique polynomial of degree zero. 
A function $f \in \ltworandom$ exhibits the polynomial chaos (PC) expansion 
\begin{equation} \label{eq:pc-f}
f(\mu) = \sum_{i=1}^{\infty} f_i \Phi_i(\mu)
\end{equation} 
with coefficients $f_i \in \R$. 
The series~(\ref{eq:pc-f}) converges in the $\ltworandom$-norm. 
In the case of time-dependent functions $f(t,\mu)$ for $t \in I$ with $I \subseteq \R$, 
the PC expansion is used pointwise for each $t \in I$. 
We also apply the PC expansion to vector-valued functions by considering 
each component separately. 

We employ truncated polynomial chaos expansions of the state variables 
as well as the inputs in the system~(\ref{eq:system-secondorder}) 
\begin{equation} \label{eq:pc-p-u} 
p^{(s)}(t,\mu) = \sum_{i=1}^s p_i(t) \Phi_i(\mu), \qquad
u^{(s)}(t,\mu) = \sum_{i=1}^s u_i(t) \Phi_i(\mu) . 
\end{equation} 
The dimensionalities of the coefficients are 
$p_i \in \R^n$ and $u_i \in \R^{n_{\rm in}}$. 
The number of basis polynomials up to a total degree~$d$ is 
$s = \frac{(q+d)!}{q!d!}$, see~\cite[p.~65]{xiu-book}.

\begin{remark} \label{remark:input} \em
If the input~$u$ does not depend on the random parameters~$\mu$, 
then still the truncated expansion~(\ref{eq:pc-p-u}) has to be used 
to achieve the same number of inputs and outputs in the final 
stochastic Galerkin system. 
Later $u_i \equiv 0$ for $i>1$ can be chosen to apply
a parameter-independent input. 
\end{remark}

\subsection{Stochastic Galerkin Method}
We introduce the following notation.

\begin{definition} \label{def:function-space}
The set of all measurable functions 
$A : \mathcal{M} \rightarrow \R^{n \times m}$, $A = (a_{k\ell})$,  
such that the expected values, cf.~(\ref{eq:expected-value}),
$$ \hat{a}_{ijk\ell} = \mathbb{E} \left[ a_{k\ell} \Phi_i \Phi_j \right] = 
\int_{\mathcal{M}} a_{k\ell}(\mu) \Phi_i(\mu) \Phi_j(\mu) \rho(\mu) \; 
\mathrm{d}\mu $$
are finite for all $i,j \in \N$ and $1 \le k \le n$, $1 \le \ell \le m$, 
is denoted by $\mathcal{F}^{n,m}$. 
The associated stochastic Galerkin projection of $A \in \mathcal{F}^{n,m}$ 
with $s$ modes reads as 
\begin{equation} \label{eq:projection-matrix} 
\hat{A} = \mathcal{G}_s(A) ,
\end{equation}
where $\hat{A} = (\hat{A}_{ij})_{i,j = 1,\ldots,s}$ comprises the submatrices
$\hat{A}_{ij} = (\hat{a}_{ijk\ell}) \in \R^{n \times m}$.
\end{definition}

We assume $M,D,K \in \mathcal{F}^{n,n}$, $B' \in \mathcal{F}^{n,n_{\rm in}}$, and 
$F,G \in \mathcal{F}^{n_{\rm out},n}$ in the system~(\ref{eq:system-secondorder}). 
The stochastic Galerkin method yields the larger system of second order
\begin{align} \label{eq:system-galerkin}
\begin{split}
\hat{M} \ddot{\hat{p}} + \hat{D} \dot{\hat{p}} + \hat{K} \hat{p} & = 
\hat{B}' \hat{u} \\[1ex] 
\hat{F} \hat{p} + \hat{G} \dot{\hat{p}} & = \hat{y}
\end{split}
\end{align}
with $\hat{p} = (\hat{p}_1^\top,\ldots,\hat{p}_s^\top)^\top \in \R^{ns}$ and 
$\hat{u} = (u_1^\top,\ldots,u_s^\top)^\top \in \R^{n_{\rm in}s}$. 
The involved matrices are $\hat{M} = \mathcal{G}_s(M)$, etc. 
using~(\ref{eq:projection-matrix}).
The coefficients in~$\hat{u}$ are the same as in~(\ref{eq:pc-p-u}). 
However, a solution~$\hat{p}$ of~(\ref{eq:system-galerkin}) yields just 
an approximation of the coefficients in~(\ref{eq:pc-p-u}), 
since the Galerkin method introduces an additional error. 
The approximation of the random-dependent state variables reads as 
\begin{equation} \label{eq:galerkin-statevariables} 
\tilde{p}^{(s)}(t,\mu) = \sum_{i=1}^s \hat{p}_i(t) \Phi_i(\mu) 
\end{equation}
for $t \ge 0$ and $\mu \in \mathcal{M}$. 

\begin{lemma} \label{lemma:definite}
Let $A \in \mathcal{F}^{n,n}$. 
If $A(\mu) \in \mathcal{S}_{\succ}^n$ for almost all $\mu \in \mathcal{M}$, 
then the stochastic Galerkin projection~(\ref{eq:projection-matrix}) 
satisfies $\hat{A} \in \mathcal{S}_{\succ}^{ns}$. 
Likewise, the property $A(\mu) \in \mathcal{S}_{\succeq}^n$ for almost all $\mu \in \mathcal{M}$ 
implies $\hat{A} \in \mathcal{S}_{\succeq}^{ns}$. 
If $A(\mu)$ is skew-symmetric for almost all~$\mu \in \mathcal{M}$, 
then $\hat{A}$ is also skew-symmetric.
\end{lemma}

The proof of the positive definite case is included in~\cite{pulch-jmi}. 
The positive semi-definite case as well as the skew-symmetric case 
can be shown by a similar argument.
 
Now we arrange a pH stochastic Galerkin system 
as in the deterministic case.

\begin{theorem} \label{thm:sgal}
Given the stochastic Galerkin system~(\ref{eq:system-galerkin}), 
an associated pH system~(\ref{eq:ph}) can be obtained 
for $\hat{x} = (\dot{\hat{p}},\hat{p})^\top$ using the matrices
$$ \hat{E} = 
\begin{pmatrix} \hat{M} & 0 \\ 0 & I_{ns} \\ \end{pmatrix} , \qquad 
\hat{J} = 
\begin{pmatrix} 0 & -I_{ns} \\ I_{ns} & 0 \\ \end{pmatrix} , \qquad 
\hat{R} = 
\begin{pmatrix} \hat{D} & 0 \\ 0 & 0 \\ \end{pmatrix} $$
$$ \hat{Q} = 
\begin{pmatrix} I_{ns} & 0 \\ 0 & \hat{K} \\ \end{pmatrix} , \qquad
\hat{B} = 
\begin{pmatrix} \hat{B}' \\ 0 \\ \end{pmatrix} , $$
and $\hat{P}=0$, $\hat{S}=0$, $\hat{N}=0$.
\end{theorem}

Proof: \nopagebreak 

It holds that $\hat{M},\hat{K} \in \mathcal{S}_{\succ}^n$ and 
$\hat{D} \in \mathcal{S}_{\succeq}^n$ due to Lemma~\ref{lemma:definite}. 
Obviously, the matrices~$\hat{J}$ and $\hat{N}$ are skew-symmetric. 
Furthermore, it follows that
$$ \hat{E}^\top \hat{Q} = 
\begin{pmatrix} \hat{M} & 0 \\ 0 & I_{ns} \\ \end{pmatrix}^\top 
\begin{pmatrix} I_{ns} & 0 \\ 0 & \hat{K} \\ \end{pmatrix} = 
\begin{pmatrix} \hat{M} & 0 \\ 0 & \hat{K} \\ \end{pmatrix} 
\in \mathcal{S}_{\succ}^{2ns} . $$
The matrix~(\ref{eq:matrix-w}) becomes
$$ \hat{W} = 
\begin{pmatrix} 
\hat{Q}^\top \hat{R} \hat{Q} & \hat{Q}^\top \hat{P} \\ 
\hat{P}^\top \hat{Q} & \hat{S} \\ 
\end{pmatrix} =
\begin{pmatrix} 
\hat{D} & 0 & 0 \\ 0 & 0 & 0 \\ 0 & 0 & 0 \\
\end{pmatrix} 
\in \mathcal{S}_{\succeq}^{2ns + n_{\rm out}} . $$
Hence the requirements of a pH system are satisfied. 
\hfill $\Box$

\medskip

The pH system from Theorem~\ref{thm:sgal} exhibits the output
$$ \hat{y} = \hat{B}^\top \hat{Q} \hat{x} = 
\begin{pmatrix} \hat{B}'^\top & 0 \\ 0 & 0 \\ \end{pmatrix}
\begin{pmatrix} \dot{\hat{p}} \\ \hat{p} \\ \end{pmatrix} 
= \hat{B}'^\top \dot{\hat{p}} . $$ 
Thus the outputs include only first-order derivatives of $\hat{p}$. 

\subsection{Hamiltonian Function of Stochastic Galerkin System}
The deterministic Hamiltonian function of the pH stochastic Galerkin system 
becomes, cf.~(\ref{eq:hamiltonian}), 
\begin{equation} \label{eq:hamiltonian-galerkin} 
\hat{H}(\hat{x}(t)) = 
\tfrac{1}{2} \hat{x}(t)^\top \hat{E}^\top \hat{Q} \hat{x}(t) = 
\tfrac{1}{2} \left( \dot{\hat{p}}(t)^\top \hat{M} \dot{\hat{p}}(t) + 
\hat{p}(t)^\top \hat{K} \hat{p}(t) \right) . 
\end{equation}
The following theoretical result demonstrates that this Hamiltonian 
can be interpreted as an expected value.

\begin{theorem} \label{thm:hamiltonian-expectation}
The Hamiltonian function~(\ref{eq:hamiltonian-galerkin}) of the 
stochastic Galerkin system coincides with the expected value
\begin{equation} \label{eq:hamiltonian-expectation} 
\hat{H}(\hat{x}(t)) = \mathbb{E} 
\left[ \tfrac{1}{2} \left( 
\dot{\tilde{p}}^{(s)}(t,\cdot) M(\cdot) \dot{\tilde{p}}^{(s)}(t,\cdot) 
+ \tilde{p}^{(s)}(t,\cdot) K(\cdot) \tilde{p}^{(s)}(t,\cdot) 
\right) \right] 
\end{equation}
for each $t \ge 0$ 
including the approximation~(\ref{eq:galerkin-statevariables}).
\end{theorem}

Proof: \nopagebreak

We compute 
\begin{align*} 
\hat{p}(t)^\top \hat{K} \hat{p}(t) & = 
\sum_{i,j=1}^s \hat{p}_i(t)^\top \hat{K}_{ij} \hat{p}_j(t) \\[0.1ex]
& = \sum_{i,j=1}^s \hat{p}_i(t)^\top 
\bigg( \int_{\mathcal{M}} K(\mu) \Phi_i(\mu) \Phi_j(\mu) \rho(\mu) \; 
\mathrm{d}\mu \bigg) \hat{p}_j(t) \\[0.1ex]
& = \int_{\mathcal{M}} 
\bigg( \sum_{i=1}^s \hat{p}_i(t) \Phi_i(\mu) \bigg)^\top 
K(\mu) \bigg( \sum_{j=1}^s \hat{p}_j(t) \Phi_j(\mu) \bigg) 
\, \rho(\mu) \; \mathrm{d}\mu \\[0.1ex]
& = 
\int_{\mathcal{M}} \tilde{p}^{(s)}(t,\mu)^\top K(\mu) \tilde{p}^{(s)}(t,\mu) 
\, \rho(\mu) \; \mathrm{d}\mu 
\; = \; \mathbb{E} 
\left[ \tilde{p}^{(s)}(t,\cdot) K(\cdot) \tilde{p}^{(s)}(t,\cdot) \right] ,
\end{align*}
where the probabilistic integration is performed separately 
in each component.
Likewise, these calculations apply to the term with $\hat{M}$. 
Basic computations generate the 
formula~(\ref{eq:hamiltonian-expectation}). 
\hfill $\Box$

\medskip

\begin{remark} 
The statement of Theorem~\ref{thm:hamiltonian-expectation} can be 
generalised to linear pH systems. 
It holds that 
$$ \hat{H}(\hat{x}(t)) = 
\mathbb{E} \left[ H(\tilde{x}^{(s)}(t,\cdot),\cdot) \right] $$
for each $t \ge 0$ under some assumptions.
Therein, $\hat{x}$ represents the solution of a stochastic Galerkin system 
and $\tilde{x}^{(s)}$ denotes the associated approximation of the 
random state variables in the original system.
\end{remark}

We assume that the stochastic Galerkin method is convergent 
for both the state variables and their first-order derivatives. 
In this case, the state variables $p(t,\mu)$ of the original dynamical 
system~(\ref{eq:system-secondorder}) and the 
approximation~(\ref{eq:galerkin-statevariables}) satisfy 
\begin{align*} 
\lim_{s \rightarrow \infty} 
\left\| p_j(t,\cdot) - \tilde{p}_j^{(s)}(t,\cdot) 
\right\|_{\ltworandom} & = 0 \\[1ex]
\lim_{s \rightarrow \infty} 
\left\| \dot{p}_j(t,\cdot) - \dot{\tilde{p}}_j^{(s)}(t,\cdot) 
\right\|_{\ltworandom} & = 0
\end{align*}
pointwise for each $t \ge 0$ and $j=1,\ldots,n$. 
The speed of convergence depends on the smoothness of the 
functions with respect to the random parameters~$\mu$, see~\cite{xiu-book}.  
Typically, the solutions of~(\ref{eq:system-secondorder}) 
inherit the differentiability of the matrices $M(\mu),D(\mu),K(\mu)$. 
In view of~(\ref{eq:hamiltonian-deterministic}),
Theorem~\ref{thm:hamiltonian-expectation} implies the convergence
\begin{equation} \label{eq:hamiltonian-convergence} 
\lim_{s \rightarrow \infty} 
\hat{H}(\hat{x}(t)) = \mathbb{E} 
\left[ H(x(t,\cdot),\cdot) \right] 
\end{equation}
for $t \ge 0$ along the transient state variables of the systems. 
We note that the matrices $\hat{E},\hat{Q}$ as well as the state variables $\hat{x}$, 
which are included in the Hamiltonian function~$\hat{H}$, 
depend on the number~$s$ of stochastic modes.

The general formulas~(\ref{eq:hamiltonian}) and~(\ref{eq:hamiltonian-bound}) 
yield
$$ \hat{H}(\hat{x}(t_2)) - \hat{H}(\hat{x}(t_1)) = 
\int_{t_1}^{t_2} \hat{y}(t)^\top \hat{u}(t) - 
\dot{\hat{p}}(t)^\top \hat{D} \dot{\hat{p}}(t) \; \mathrm{d}t \le 
\int_{t_1}^{t_2} \hat{y}(t)^\top \hat{u}(t) \; \mathrm{d}t $$
for $t_2 \ge t_1 \ge 0$. 
Now we consider a deterministic input by $\hat{u}(t) = (u_1(t),0,\ldots,0)^\top$ 
with $u_1 \in \R^{n_{\rm in}}$, 
i.e., the input does not depend on the random variables, 
see Remark~\ref{remark:input}.
It follows that 
$$ \hat{H}(\hat{x}(t_2)) - \hat{H}(\hat{x}(t_1)) \le 
\int_{t_1}^{t_2} \hat{y}_1(t)^\top u_1(t) \; \mathrm{d}t , $$
where $\hat{y}_1$ represents an approximation of the expected value 
of the output in the original dynamical 
system~(\ref{eq:system-secondorder}) with $F=0$ and $G=B'^\top$. 

For simplicity, we consider an SISO system~(\ref{eq:system-secondorder}) 
with input~$u$ and output~$y$ in the following. 
The approximation reads as 
\begin{align*} 
\int_{t_1}^{t_2} \hat{y}_1(t) u_1(t) \; \mathrm{d}t & \approx
\int_{t_1}^{t_2} \mathbb{E} [ y(t,\cdot) ] \, u(t) \; \mathrm{d}t = 
\int_{t_1}^{t_2} \int_{\mathcal{M}} y(t,\mu) u(t) \rho(\mu) \; 
\mathrm{d}\mu \mathrm{d}t \\[1ex]
& = \int_{\mathcal{M}} \int_{t_1}^{t_2} y(t,\mu) u(t) \rho(\mu) \; 
\mathrm{d}t \mathrm{d}\mu 
= \mathbb{E} \bigg[ \int_{t_1}^{t_2} y(t,\cdot) u(t) \; \mathrm{d}t \bigg] ,
\end{align*}
where the two integrations can be interchanged due to Fubini's theorem. 
Therein, we obtain the expected value of the upper bound for the 
Hamiltonian function~(\ref{eq:hamiltonian-deterministic}) 
belonging to the system~(\ref{eq:system-secondorder}).


\section{Model Order Reduction}
\label{sec:mor} 
We investigate structure-preserving MOR of second-order ODEs, 
see~\cite{chu-tsai-lai,reis-stykel,salimbahrami-lohmann}. 
The stochastic Galerkin system~(\ref{eq:system-galerkin}) represents 
our full-order model (FOM), since the dimensionality is high 
in the case of many random variables. 
A reduced-order model (ROM) reads as
\begin{align} \label{eq:system-mor}
\begin{split}
\bar{M} \ddot{\bar{p}} + \bar{D} \dot{\bar{p}} + \bar{K} \bar{p} & = 
\bar{B}' \hat{u} \\[1ex]
\bar{F} \bar{p} + \bar{G} \dot{\bar{p}} & = \bar{y} 
\end{split}
\end{align}
with state variables $\bar{p} = (\bar{p}_1^\top,\ldots,\bar{p}_r^\top)^\top$ 
and the same input $\hat{u} = (u_1^\top,\ldots,u_s^\top)^\top$ as 
in~(\ref{eq:system-galerkin}).  
Now the dimensionalities are $\bar{p}_i \in \R$ for each~$i$ and $u_i \in \R^{n_{\rm in}}$. 

In projection-based MOR, two projection matrices $V,W \in \R^{ns \times r}$ are 
computed, where $V$ typically is an orthogonal matrix ($V^\top V = I_r$). 
Often the biorthogonality condition $W^\top V = I_r$ is added. 
In the case of $W = V$, the reduction is called a Galerkin-type MOR method. 
In the case of $W \neq V$, the reduction is referred to as a 
Petrov-Galerkin MOR method. 

\begin{definition} \label{def:reduced-matrix}
Let $A \in \R^{n \times n}$ be a constant square matrix 
and $V,W \in \R^{n \times r}$ be two projection matrices. 
The projection operator for the matrix~$A$ is given by 
\begin{equation} \label{eq:reduction-operator} 
\bar{A} = \mathcal{P}_r(A) = W^\top A V 
\end{equation}
with the result $\bar{A} \in \R^{r \times r}$. 
\end{definition}

Concerning the ROM~(\ref{eq:system-mor}), it holds that 
$\bar{A} = \mathcal{P}_r(A)$ for $A \in \{ \hat{M} , \hat{D} , \hat{K} \}$, 
$\bar{B}' = W^\top \hat{B}'$, and
$\bar{A} = A V$ for $A \in \{ \hat{F} , \hat{G} \}$. 
The combination of the operators in 
Definition~\ref{def:function-space} and 
Definition~\ref{def:reduced-matrix} yields
$$ \bar{A} = \mathcal{P}_r(\hat{A})
= \mathcal{P}_r(\mathcal{G}_s(A)) $$
for matrices like $A \in \{M,D,K\}$ 
in the second-order system~(\ref{eq:system-secondorder}). 
The numbers $r$ and $s$ are independent of each other 
(just $r < ns$ is required). 
Now the structure of the original matrices should be preserved. 

A symmetric matrix becomes an unsymmetric matrix in a 
Petrov-Galerkin projection~(\ref{eq:reduction-operator}). 
Alternatively, both symmetry and (semi-)definiteness of matrices 
is preserved in an MOR with Galerkin projection. 
A semi-definite matrix may become definite in this projection. 
It follows that $\bar{M},\bar{K} \in \mathcal{S}_{\succ}^r$, and 
$\bar{D} \in \mathcal{S}_{\succeq}^r$. 
Thus a Galerkin-type MOR of the second-order system preserves 
the structure. 
Moreover, the ROM~(\ref{eq:system-mor}) is asymptotically stable 
in the case of $\bar{D} \in \mathcal{S}_{\succ}^r$, 
which is guaranteed for $\hat{D} \in \mathcal{S}_{\succ}^{ns}$. 

Typical MOR schemes for second-order systems are Krylov subspace methods, 
see~\cite{salimbahrami-lohmann}, and balanced truncation, 
see~\cite{reis-stykel}. 
An elementary Krylov subspace method is the Arnoldi technique, 
which yields a Galerkin-type MOR. 
Hence this method preserves the structure of our second-order 
stochastic Galerkin system.


\section{Test Example}
\label{sec:example}
We investigate the properties, 
which are reported in the previous sections, 
using a test example of a second-order linear dynamical system.

\subsection{Modelling}
\begin{figure}
  \begin{center}
  \includegraphics[width=5cm]{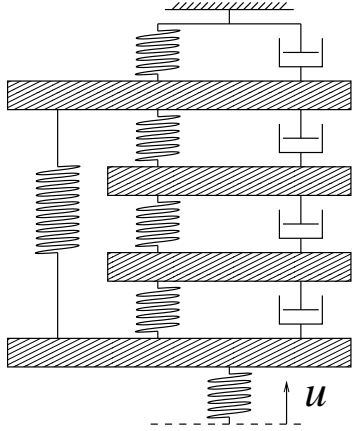}
  \end{center}
  \caption{Mass-spring-damper system.}
\label{fig:msd}
\end{figure}

\begin{figure}
  \begin{center}
  \includegraphics[width=6.5cm]{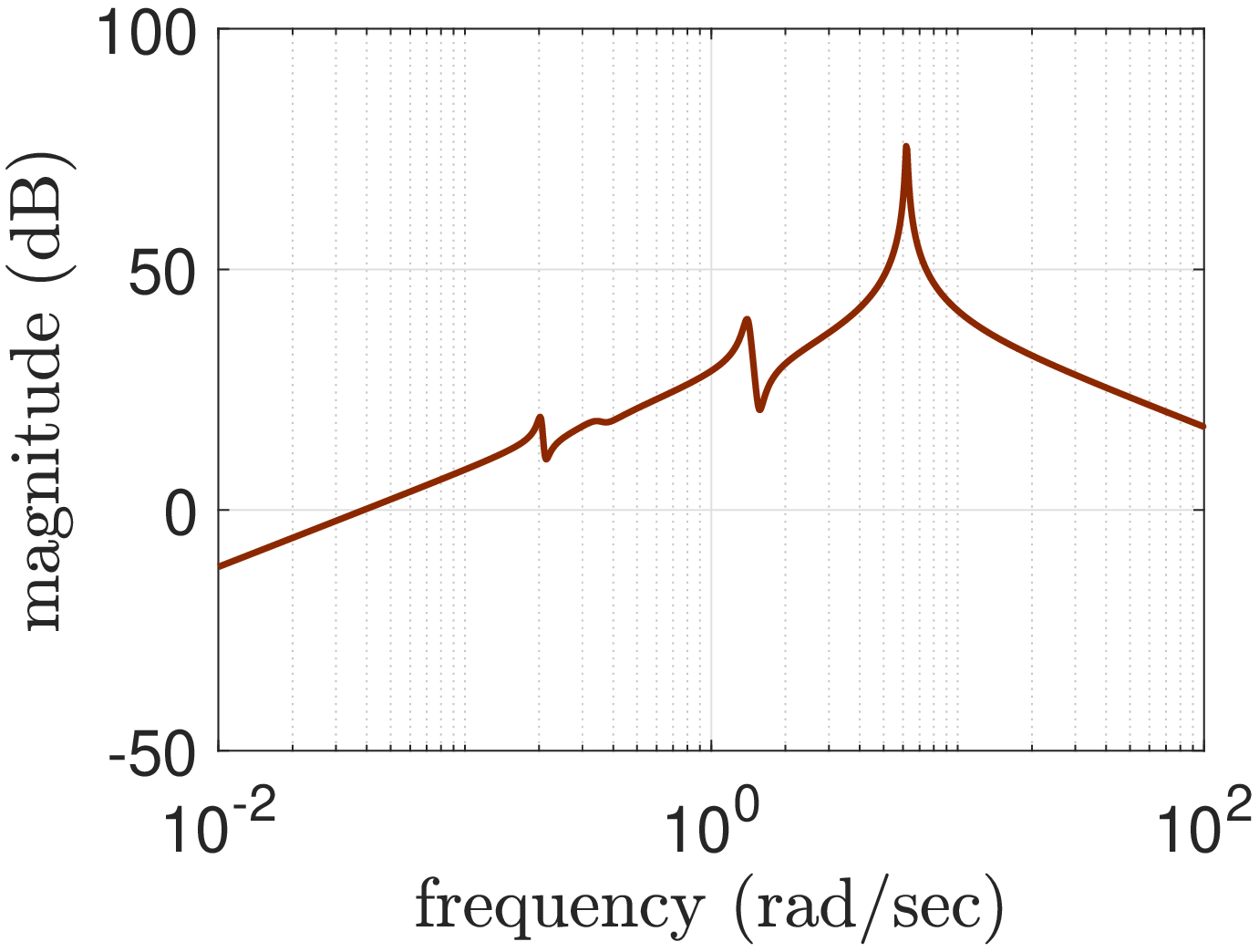}
  \hspace{5mm}
  \includegraphics[width=6.5cm]{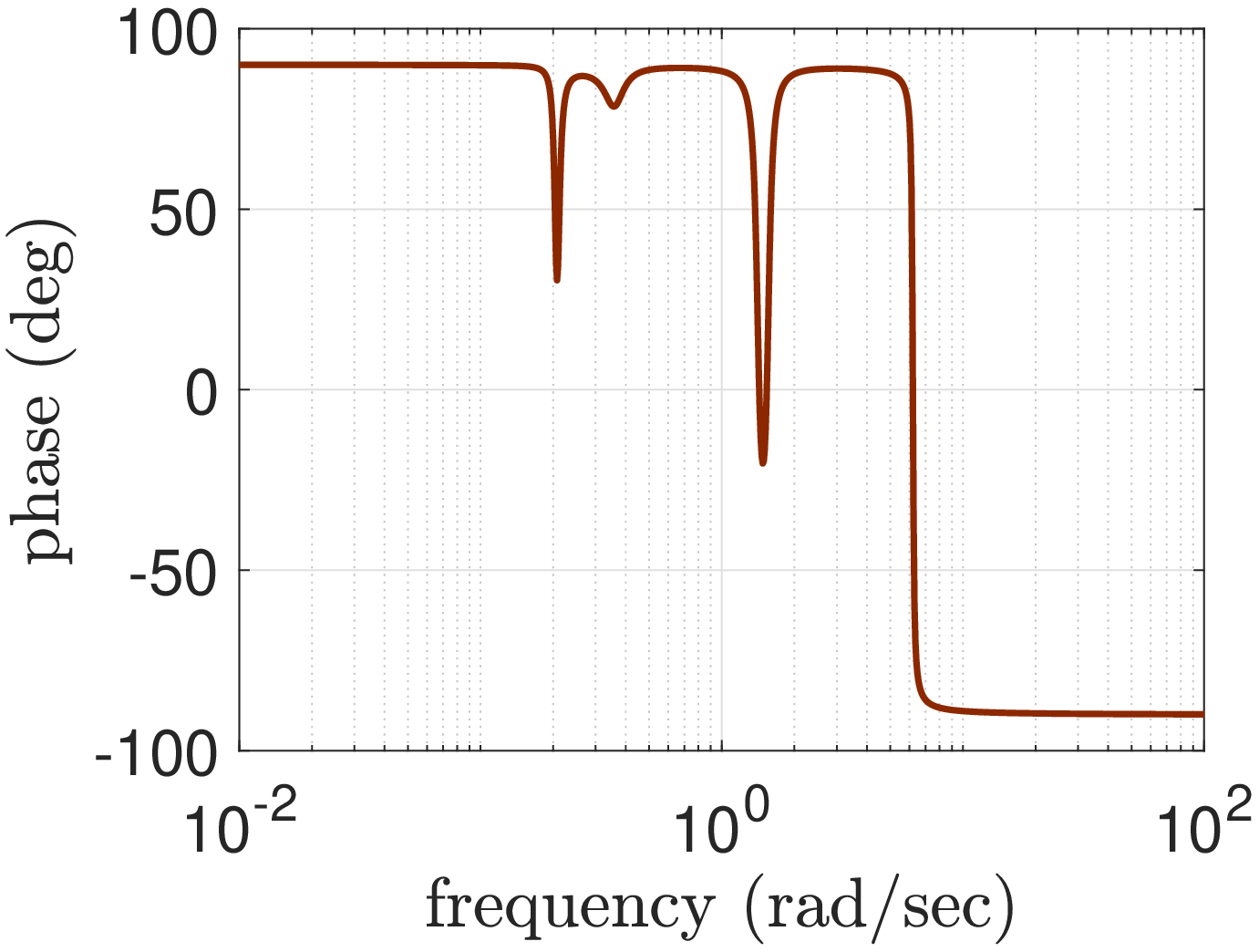}
  \end{center}
  \caption{Bode plot of mass-spring-damper system for a 
           constant choice of the physical parameters.}
\label{fig:bode}
\end{figure}

We study a mass-spring-damper system from~\cite{lohmann-eid}, 
depicted in Figure~\ref{fig:msd}. 
There are $14$ physical parameters in this configuration: 
four masses, four damping constants and six spring constants.  
A mathematical modelling yields a system of four second-order ODEs 
in the form~(\ref{eq:system-secondorder}). 
The mass matrix, the damping matrix, and the stiffness matrix 
are symmetric and positive definite provided that 
all physical parameters are positive. 
Thus the linear dynamical system is asymptotically stable.  
The single input~$u$ is an excitation at the bottom spring. 
This input implies that the single output~$y$ is the velocity 
of the bottom mass in the associated pH system. 
Hence an SISO system is given.
Figure~\ref{fig:bode} shows the Bode plot of this linear dynamical system 
for a deterministic selection of the parameters. 

In a stochastic modelling, we replace the $q=14$ physical parameters by 
independent random variables with uniform distributions varying 
10\% around their mean values. 
The PC expansions incorporate multivariate basis polynomials, 
which are products of the univariate Legendre polynomials. 
We apply truncated PC expansions including all polynomials up to 
total degree $d=3$. 
It follows that the number of basis polynomials is $s=680$. 
We arrange the stochastic Galerkin system~(\ref{eq:system-galerkin}), 
which exhibits dimension $ns = 2720$. 
The number of inputs in~$\hat{u}$ as well as the 
number of outputs in~$\hat{y}$ is equal to~$s$. 

\subsection{Transient Simulation}
\label{sec:transient-simulation}
In this section, we examine solutions of initial value problems (IVPs) 
in a time interval $[0,100]$, 
where the initial conditions are always chosen to be zero. 
Second-order ODE systems are converted into equivalent 
first-order ODE systems, such that numerical methods for first-order systems 
are employed. 

\begin{figure}
  \begin{center}
  \includegraphics[width=6.5cm]{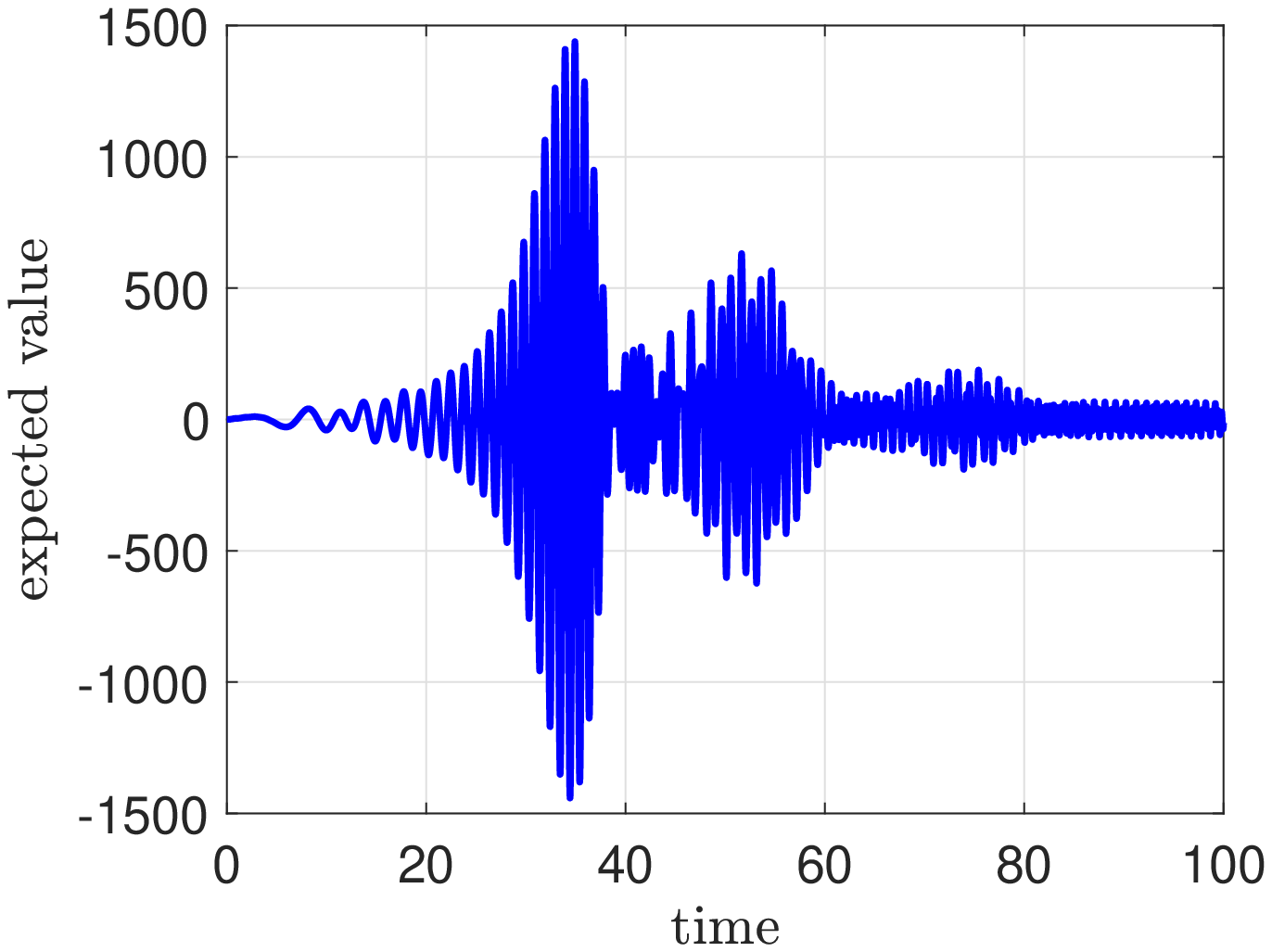}
  \hspace{5mm}
  \includegraphics[width=6.5cm]{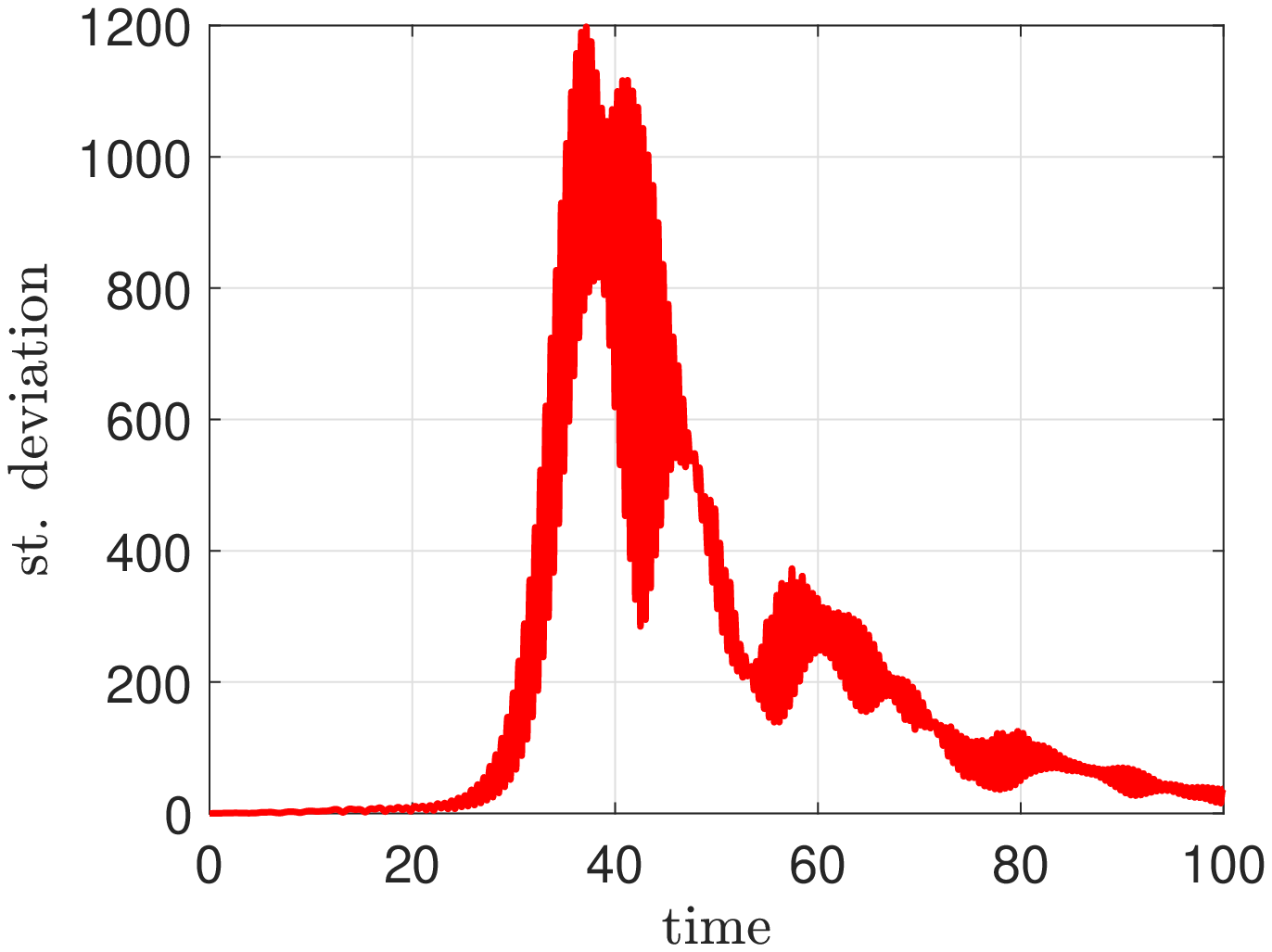}
  \end{center}
  \caption{Approximations of expected value (left) and 
           standard deviation (right) for random QoI 
           in mass-spring-damper system.}
\label{fig:galerkin-statistics}
\end{figure}

We solve an IVP of the stochastic Galerkin system~(\ref{eq:system-galerkin}). 
A deterministic input $\hat{u} = (u_1,0,\ldots,0)^\top$ is supplied with the signal  
$$ u_1(t) = \sin \left( \tfrac{1}{10} t^2 \right) , $$
which can be interpreted as a harmonic oscillation with increasing frequency. 
A Runge-Kutta method of order 4(5) yields a numerical solution, 
where a step size control is applied based on a local error control 
with relative tolerance $\varepsilon_{\rm rel} = 10^{-4}$ and 
absolute tolerance $\varepsilon_{\rm abs} = 10^{-6}$. 
The numerical solution generates an approximation of the 
multiple outputs~$\hat{y}$. 
Figure~\ref{fig:galerkin-statistics} demonstrates the expected value 
as well as the standard deviation associated to the QoI of the random ODEs, 
which are obtained by the transient solution of the stochastic Galerkin 
system. 
We observe a maximum resonance effect around $t \approx 35$ 
in the expected value, where also a maximum standard deviation appears. 

\begin{figure}
  \begin{center}
  \includegraphics[width=6.5cm]{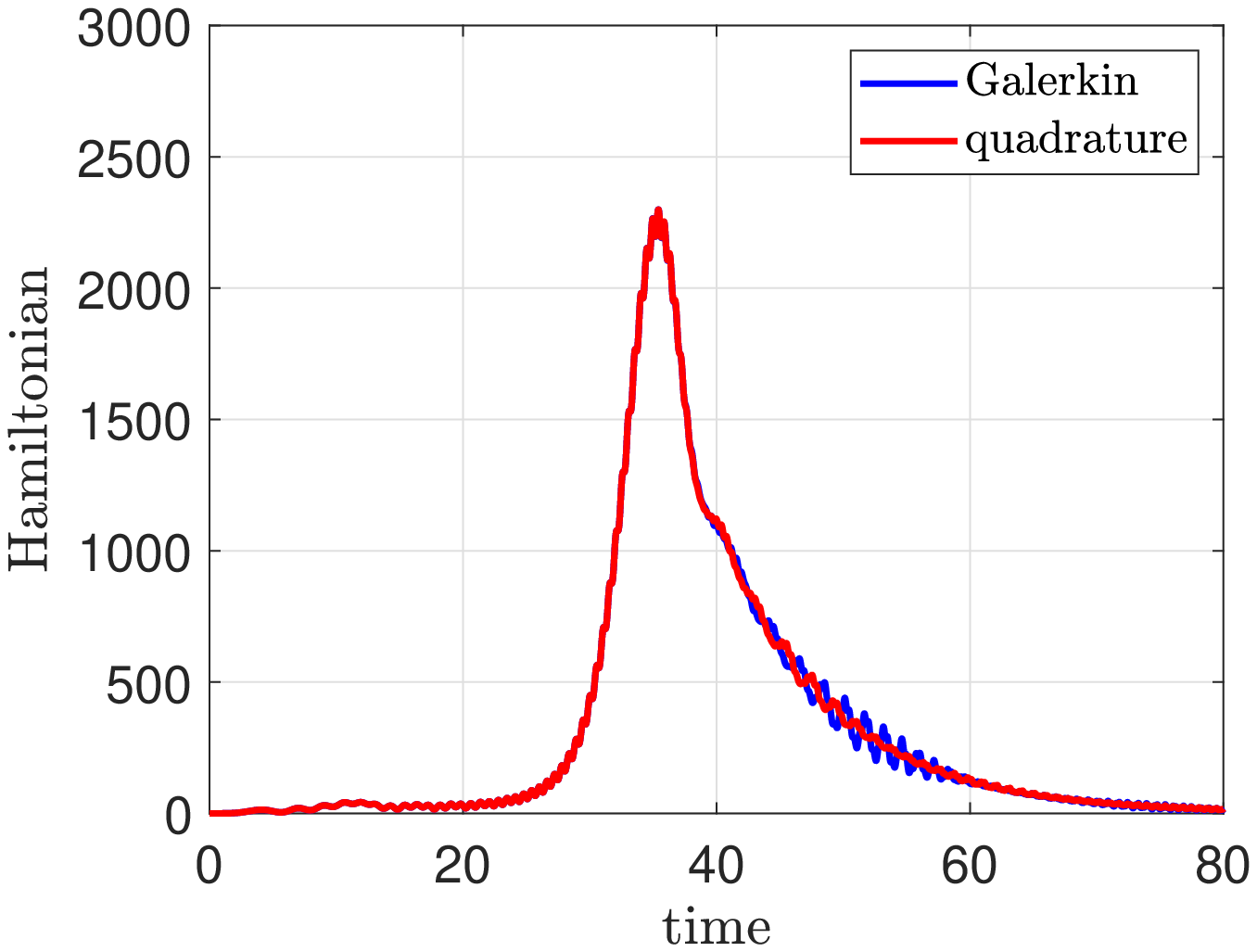}
  \hspace{5mm}
  \includegraphics[width=6.5cm]{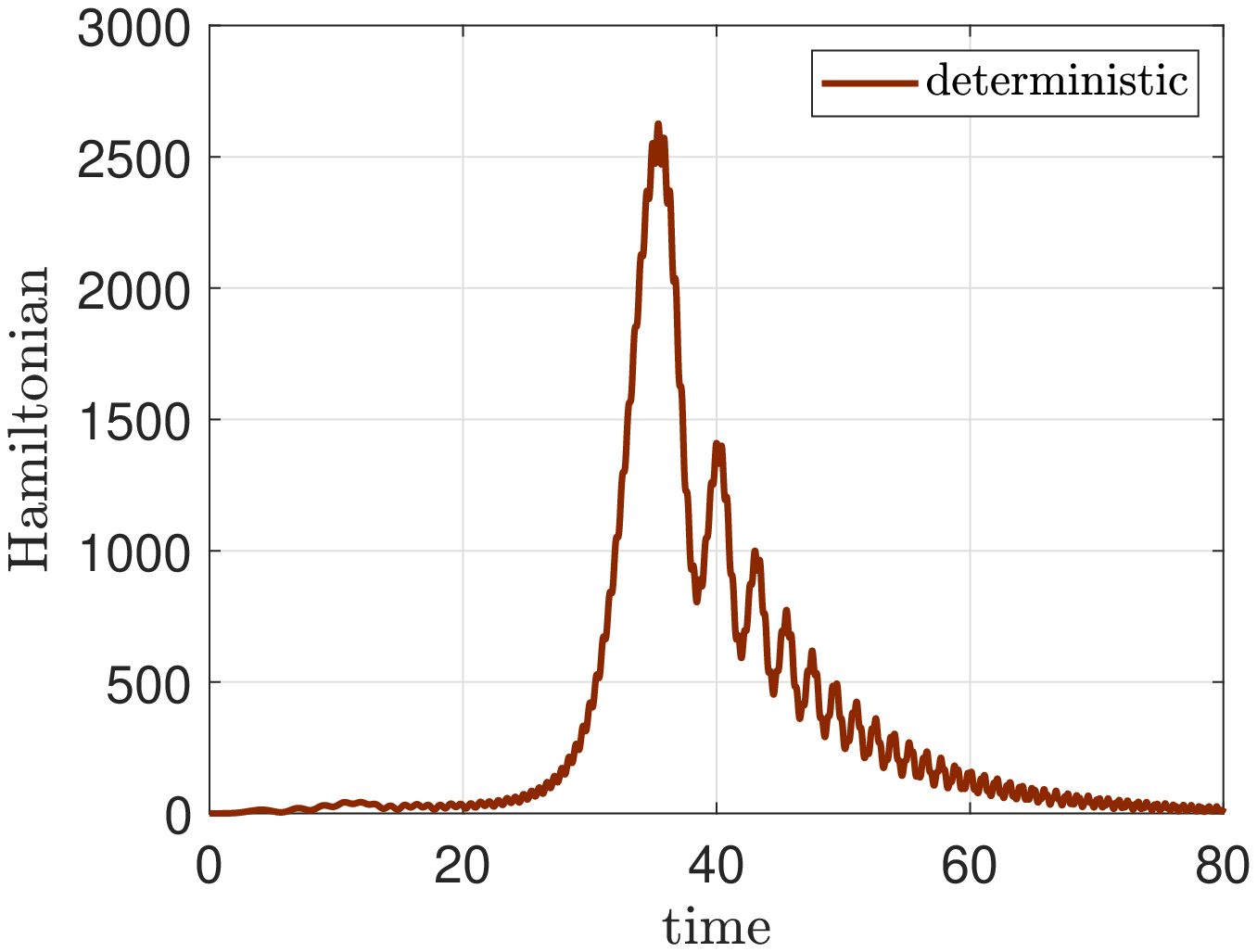}
  \end{center}
  \caption{Hamiltonian of stochastic Galerkin system as well as 
           expected value of Hamiltonians obtained by quadrature (left)
           and Hamiltonian of deterministic system (right).}
\label{fig:hamiltonian}
\end{figure}

Now we investigate the Hamiltonian functions associated to the 
transient solutions. 
Figure~\ref{fig:hamiltonian} (left) depicts the Hamiltonian belonging to 
the stochastic Galerkin system~(\ref{eq:system-galerkin}). 
We perform two additional simulations for comparison. 
On the one hand, we solve an IVP of the system~(\ref{eq:system-secondorder}) 
once using the mean value of the random variables as parameters. 
The resulting Hamiltonian function is shown in Figure~\ref{fig:hamiltonian} 
(right). 
On the other hand, we compute an approximation of the expected value 
of the Hamiltonians associated to the random 
system~(\ref{eq:system-secondorder}). 
This expected value is approximated by a quadrature using the 
Stroud-5 scheme, see~\cite{stroud}, which includes 393~nodes 
in our case of $14$~random variables. 
Hence 393 IVPs of a deterministic system~(\ref{eq:system-secondorder}) 
are solved. 
The approximation of the expected value is displayed in 
Figure~\ref{fig:hamiltonian} (left). 
We recognise a good agreement between the Hamiltonian of the 
stochastic Galerkin system and the averaged Hamiltonians from the 
quadrature method, 
which indicates the convergence~(\ref{eq:hamiltonian-convergence}). 

\subsection{Model Order Reduction}
\begin{figure}
  \begin{center}
  \includegraphics[width=6.5cm]{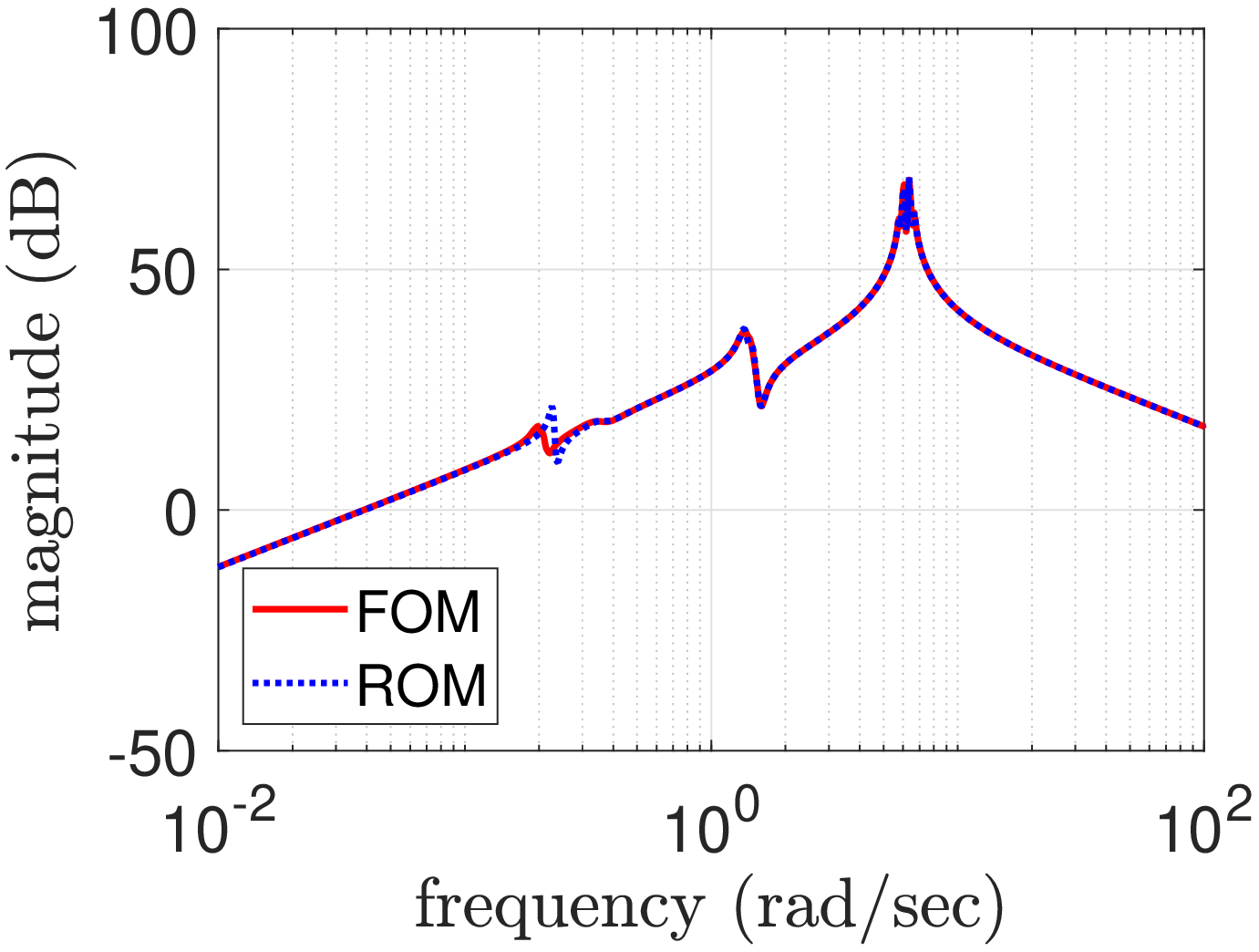}
  \hspace{5mm}
  \includegraphics[width=6.5cm]{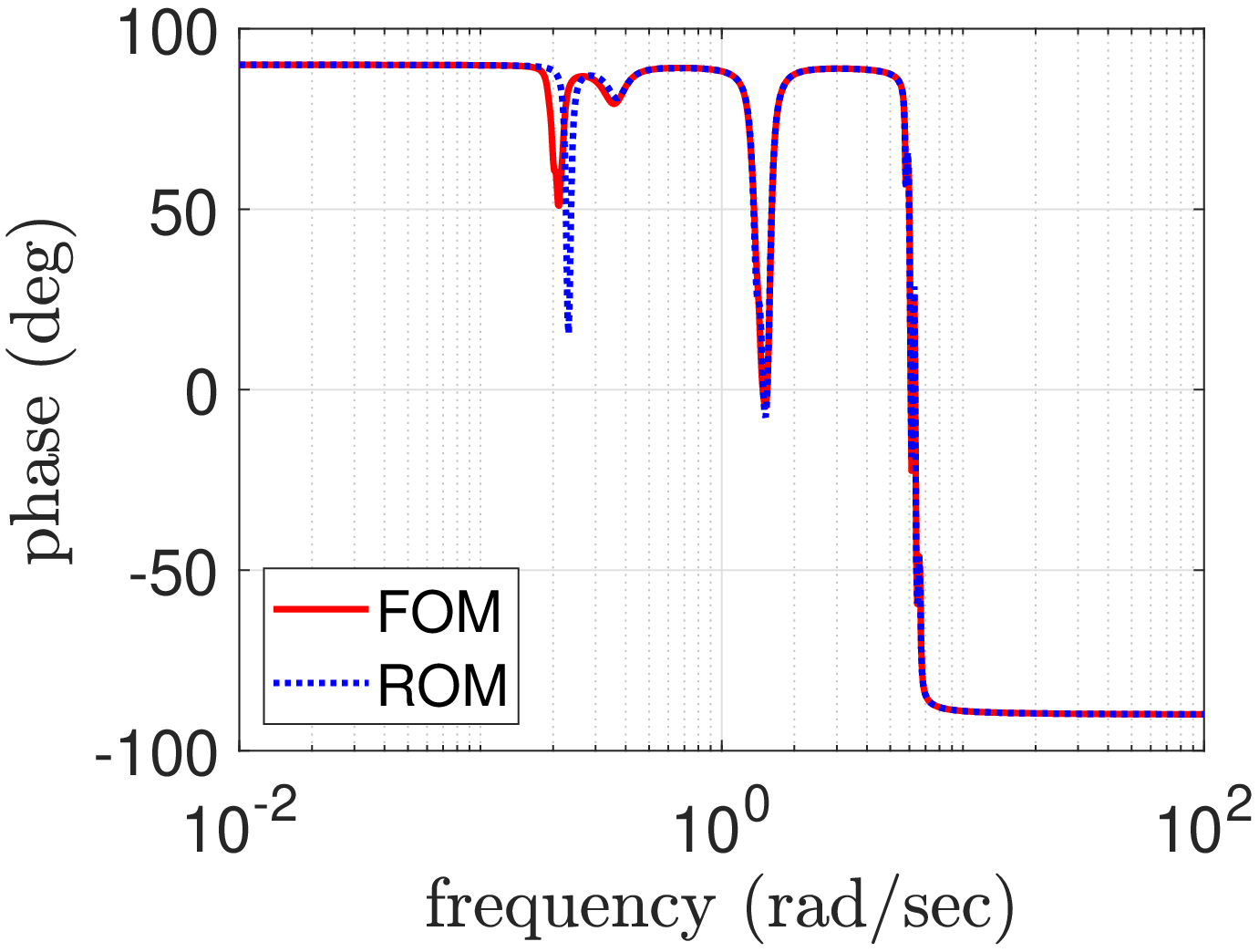}
  \end{center}
  \caption{Bode plots of FOM and an ROM of dimension~$30$.}
\label{fig:bode-reduction}
\end{figure}

We apply the second-order Arnoldi algorithm as 
in~\cite[p.~960]{chu-tsai-lai}. 
This MOR method is independent from the definition of outputs in the system. 
However, a multiple-input system requires a block Arnoldi technique, 
which becomes inefficient in the case of many inputs. 
Since we just use a single (deterministic) input, 
this problem does not occur. 
We compute a projection matrix $V$ of rank~50 by the second-order 
Arnoldi method. 
Using (only) the first $r$ columns of~$V$ yields ROMs of 
dimension $r \le 50$.
Now we restrict the output to the expected value of the QoI 
in the original system to study an SISO system.
The Bode plots of the FOM as well as the ROM of dimension $r=30$ 
are shown in Figure~\ref{fig:bode-reduction}. 
The matrices $\bar{M},\bar{D},\bar{K}$ in all 
ROMs~(\ref{eq:system-mor}) are symmetric and positive definite, 
since a structure-preserving MOR method is applied. 

\begin{figure}
  \begin{center}
  \includegraphics[width=7cm]{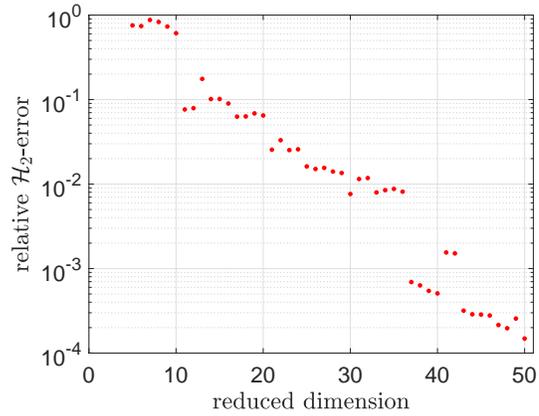}
  \end{center}
  \caption{Relative error of MOR for reduced dimensions $r=5,6,\ldots,50$ 
           using the second-order Arnoldi method.}
\label{fig:mor-error}
\end{figure}

We examine the approximation error in this MOR. 
The relative error of the associated transfer functions 
$H_{\rm FOM}$ and $H_{\rm ROM}$ reads as
$$ \frac{\| H_{\rm FOM} - H_{\rm ROM} \|_{\htwonorm}}{
\| H_{\rm FOM} \|_{\htwonorm}} $$
employing the $\htwonorm$-norm, see~\cite[p.~144]{antoulas}. 
Figure~\ref{fig:mor-error} illustrates these relative $\htwonorm$-errors 
for the ROMs of dimensions $5 \le r \le 50$. 
We recognise that the errors exponentially converge to zero for 
increasing dimensions. 
However, this convergence is not monotone. 

\begin{figure}
  \begin{center}
  \includegraphics[width=6.5cm]{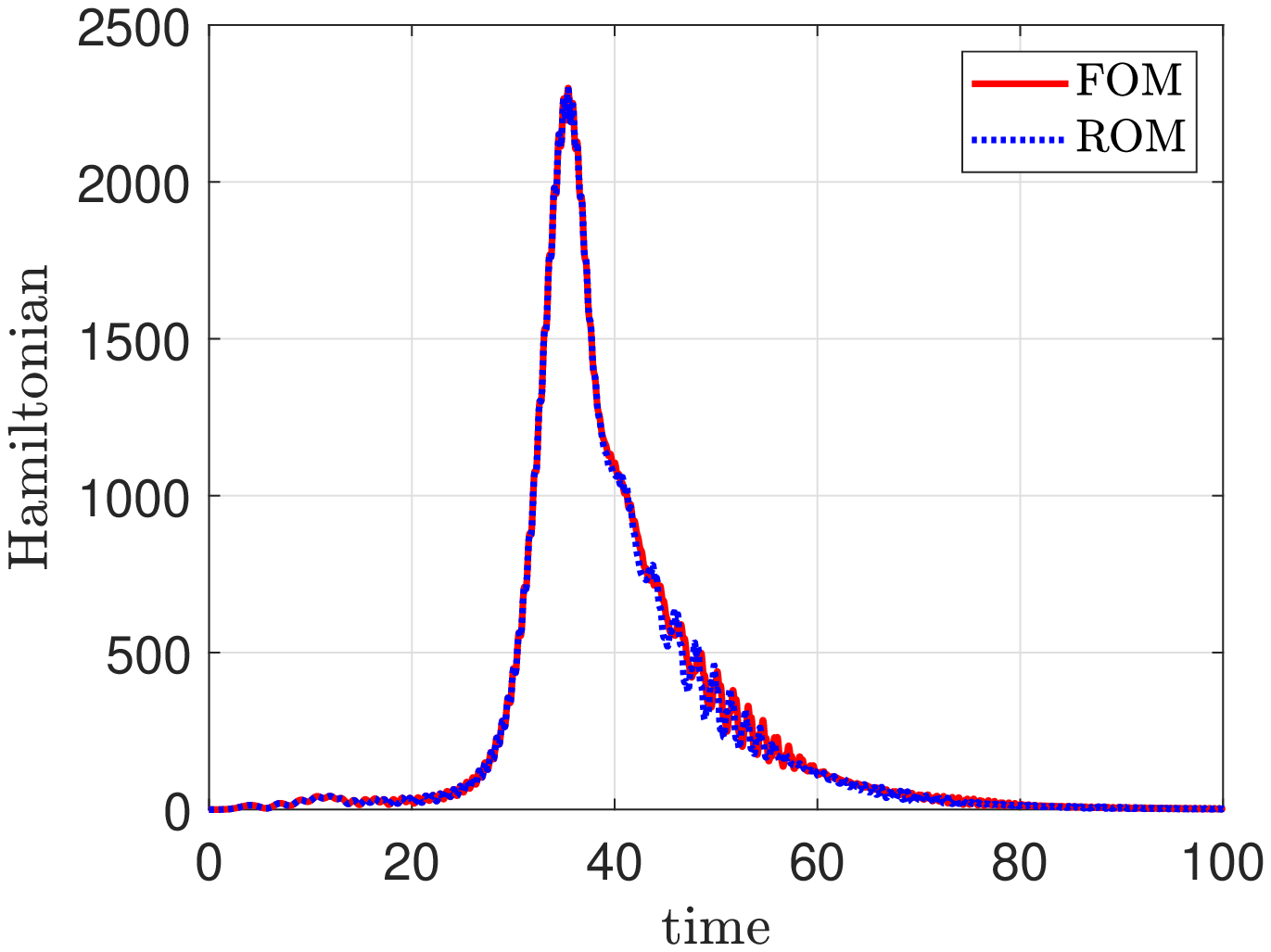}
  \hspace{5mm}
  \includegraphics[width=6.5cm]{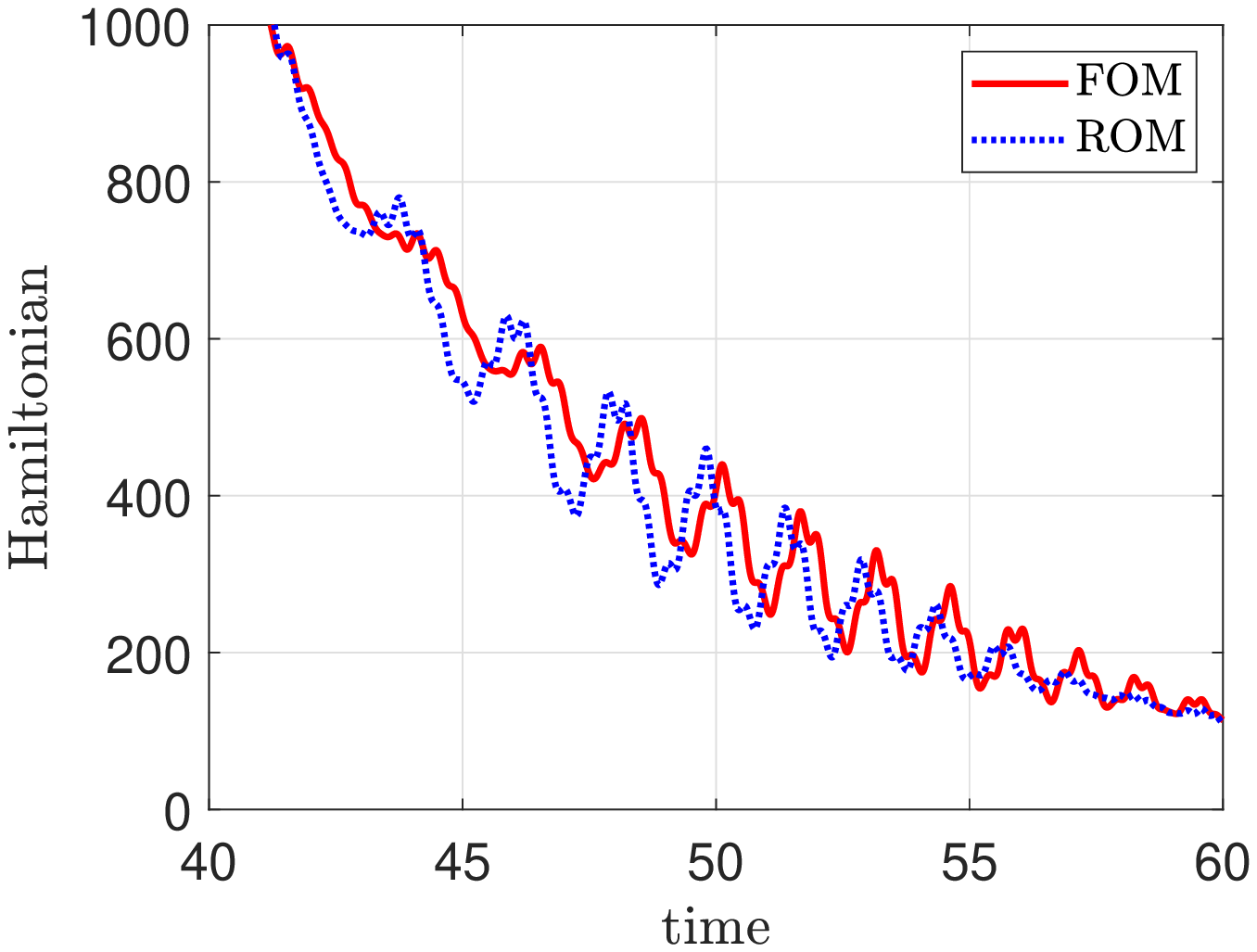}
  \end{center}
  \caption{Hamiltonian functions from transient simulation 
           of the stochastic Galerkin system (FOM) and its ROM 
           of dimension~$10$
           (left: total time interval, right: zoom).}
\label{fig:hamiltonian-mor}
\end{figure}

Furthermore, we repeat the transient simulation from 
Section~\ref{sec:transient-simulation} using the ROM 
of dimension $r=10$ instead of the full stochastic Galerkin system. 
Figure~\ref{fig:hamiltonian-mor} depicts the resulting Hamiltonian  
functions of FOM and ROM. 
If higher dimensions like $r \ge 15$ are used, then the differences 
between FOM and ROM become tiny and cannot be observed in 
a plot anymore.

\section{Conclusions}
The stochastic Galerkin method applied to a linear ODE system of 
second order yields an ODE system of the same type.  
A structure-preserving MOR of the stochastic Galerkin system generates 
smaller ODE systems of the same type again. 
We formulated a first-order pH system for each second-order system. 
The associated Hamiltonian functions were investigated. 
We showed that the Hamiltonian for the stochastic Galerkin system 
represents an approximation of the expected value of the Hamiltonian 
for the original system including random variables. 
Numerical computations confirm this theoretical result 
using a test example.  
Furthermore, we demonstrated an efficient structure-preserving MOR 
in the test example. 

\bigskip

{\bf Acknowledgement} \\
The author thanks Prof.~Dr.~Christopher~Beattie 
(Virginia Polytechnic Institute and State University, United States) 
for helpful discussions.



\begin{thebibliography}{00}

\bibitem{antoulas}
A.~Antoulas,
Approximation of Large-Scale Dynamical Systems, 
SIAM Publications, 2005.

\bibitem{bai-su}
Z.~Bai, Y.~Su, 
Dimension reduction of large-scale second-order dynamical systems 
via a second-order Arnoldi method,
SIAM J. Sci. Comput. 26:5 (2005) 1692--1709.

\bibitem{beattie-gugercin}
C.~Beattie, S.~Gugercin, 
Krylov-based model reduction of second-order systems with 
proportional damping,
in: 44th IEEE Conference on Decision and Control 
(2005) 2278--2283.

\bibitem{beattie-mehrmann}
C.~Beattie, V.~Mehrmann, H.~Xu, H.~Zwart, 
Linear port-Hamiltonian descriptor systems, 
Mathematics of Control, Signals, and Systems 30:4 (2018).

\bibitem{beattie-mehrmann-dooren}
C.~Beattie, V.~Mehrmann, P.~Van~Dooren, 
Robust port-Hamiltonian representations of passive systems, 
Automatica~100 (2019) 182-186.

\bibitem{braun}
M.~Braun, 
Differential Equations and Their Applications: 
an Introduction to Applied Mathematics, 
4th ed., 
Springer, 1993.

\bibitem{castane-selga}
R.~Casta{\~n}{\'e}-Selga, B.~Lohmann, R.~Eid, 
Stability preservation in projection-based model order reduction 
of large scale systems. 
European Journal of Control 18:2 (2012) 122-132. 

\bibitem{chu-tsai-lai}
C.-C.~Chu, H.-C.~Tsai, M.-H.~Lai, 
Structure preserving model-order reductions of MIMO second-order systems 
using Arnoldi methods, 
Math. Comput. Model. 51 (2010) 956-973.

\bibitem{hauschild}
S.-A.~Hauschild, N.~Marheineke, V.~Mehrmann, 
Model reduction techniques for linear constant coefficient 
port-Hamiltonian differential-algebraic systems, 
Control and Cybernetics 48:1 (2019).

\bibitem{inman}
D.J.~Inman, Vibration and Control, John Wiley \& Sons, Ltd, 2006.
   
\bibitem{ionescu} 
T.C.~Ionescu, A.~Astolfi, 
On moment matching with preservation of passivity and stability, 
in: 49th IEEE Conference on Decision and Control (2010), 6189-6194. 

\bibitem{jacob-zwart}
B.~Jacob, H.J.~Zwart, 
Linear Port-Hamiltonian Systems on Infinite-di\-men\-sional Spaces, 
Operator Theory: Advances and Applications, Vol.~223, 
Springer, 2012.

\bibitem{jacob-skrepek}
B.~Jacob, N.~Skrepek, 
Stability of the multidimensional wave equation in port-Hamiltonian modelling, 
in: 60th IEEE Conference on Decision and Control 
(2021) 6188-6193. 

\bibitem{lohmann-eid}
B.~Lohmann, R.~Eid, 
Efficient order reduction of parametric and nonlinear models 
by superposition of locally reduced models, 
in: G.~Roppen\-ecker, B.~Lohmann (eds.), 
Methoden und Anwendungen der Regelungs\-technik, Shaker, 2009.

\bibitem{pulch18}
R.~Pulch,
Model order reduction and low-dimensional representations for
random linear dynamical systems.
Math. Comput. Simulat. 144 (2018) 1--20.

\bibitem{pulch-jmi}
R.~Pulch, 
Stability-preserving model order reduction for 
linear stochastic Galerkin systems, 
J. Math. Ind. 9:10 (2019).

\bibitem{pulch21}
R.~Pulch, 
Frequency domain integrals for stability preservation in Galerkin-type 
projection-based model order reduction, 
Int. J. Control 94:7 (2021) 1734-1750.

\bibitem{reis-stykel} 
T.~Reis, T.~Stykel, 
Balanced truncation model reduction of second-order systems, 
Math. Comp. Mod. Dyn. Syst. 14:5 (2008), 391-406.

\bibitem{salimbahrami-lohmann}
B.~Salimbahrami, B.~Lohmann, 
Order reduction of large scale second-order systems using 
Krylov subspace methods, 
Linear Algebra and its Applications 415 (2006) 385-405.

\bibitem{schaft-jeltsema}
A.~van~der~Schaft, D.~Jeltsema, 
Port-Hamiltonian Systems Theory: An Introductory Overview, 
Foundations and Trends in Systems and Control 1:2-3 
(2014) 173-378.

\bibitem{stroud}
A.~Stroud, 
Approximate Calculation of Multiple Integrals, 
Prentice Hall, 1971.

\bibitem{sullivan}
T.J.~Sullivan, Introduction to Uncertainty Quantification,
Springer, 2015.

\bibitem{xiu-book}
D.~Xiu,  
Numerical Methods for Stochastic Computations: a Spectral Method Approach, 
Princeton University Press, 2010.

\end{thebibliography}
\end{document}